\documentclass[12 pt]{article}

\usepackage{}
\usepackage{amsmath,amsfonts,amssymb,mathrsfs,amsthm}

\theoremstyle{plain}
\newtheorem{thm}{Theorem}[section]
\newtheorem{cor}[thm]{Corollary}
\newtheorem{lem}[thm]{Lemma}
\newtheorem{prop}[thm]{Proposition}
\theoremstyle{definition}
\newtheorem{defn}{Definition}[section]
\theoremstyle{remark}
\newtheorem{rem}{\bf Remark}[section]

\begin{document}
\title{Counting of lattices containing up to $4$ reducible elements and having nullity up to $3$}
\maketitle
\begin{center}
\author{Dr. A. N. Bhavale\\hodmaths@moderncollegepune.edu.in\\Department of Mathematics, Modern College of ASC(A),\\ Shivajinagar, Pune 411005, M.S., India.}
\end{center}

\begin{abstract}
In this paper, we count all non-isomorphic lattices on $n$ elements, containing four reducible elements and having nullity three. 
This work is in respect of Birkhoff's open problem (which is NP-complete) of counting all finite lattices on $n$ elements.
\end{abstract}

MSC Classification 2020: 06A05, 06A06, 06A07.

\section{Introduction}
In 1940, Birkhoff \cite{bir} posed the following open problem.\\
{\bf Problem 1.} {\it Compute all posets or lattices on a set of $n$ elements up to isomorphism}.\\

There were attempts to solve this problem by many authors. In 2002, Brinkmann and Mckay \cite{bm} enumerated all non-isomorphic posets with $15$ and $16$ elements. The work of enumeration of all non-isomorphic (unlabeled) posets is still in progress for $n\geq 17$. In the same year, Heitzig and Reinhold \cite{hr} counted all non-isomorphic (unlabeled) lattices on up to $18$ elements. 

In 2003 Pawar and Waphare \cite{pw} counted all non-isomorphic lattices on $n$ elements and containing $n$ edges, which are precisely the lattices having nullity one.
Also in 2002 Thakare, Pawar, and Waphare \cite{tpw} counted all non-isomorphic lattices on $n$ elements, containing exactly two reducible elements.
Thakare, Pawar, and Waphare \cite{tpw} also counted all non-isomorphic lattices on $n$ elements and containing $n+1$ edges, which are precisely the lattices having nullity two.
Independently, Bhavale and Aware \cite{ba2} counted all non-isomorphic lattices on $n$ elements, having nullity two.
Recently, Bhavale and Aware \cite{ba} counted all non-isomorphic lattices on $n$ elements, containing exactly three reducible elements.
Aware and Bhavale \cite{AB} also counted all non-isomorphic lattices on $n$ elements, containing four comparable reducible elements, and having nullity three.
In this paper, we count all non-isomorphic lattices on $n$ elements, containing up to four reducible elements, and having nullity up to three.

\section{Preliminaries}
Let $\leq$ be a partial order relation on a non-empty set $P$, and let $(P,\leq)$ be a poset. Elements $x,y \in P$ are said to be $comparable$, if either $x \leq y$ or $y \leq x$. A poset is called a {\it chain} if any two elements in it are comparable.
Elements $x, y \in P$ are said to be $incomparable$, denoted by $x||y$, if $x,y$ are not comparable. An element $c\in P$ is {\it{a lower bound}} ({\it{an upper bound}}) of $a,b\in P$ if $c\leq a,c\leq b(a\leq c,b\leq c)$. A {\it{meet}} of $a,b \in P$, denoted by $a \wedge b$, is defined as the greatest lower bound of $a$ and $b$. A {\it{join}} of $a,b\in P$, denoted by $a \vee b$, is defined as the least upper bound of $a$ and $b$. A poset $L$ is a $lattice$ if $a \wedge b$ and $a \vee b$, exist in $L$, $\forall a,b\in L$. Lattices $L_1$ and $L_2$ are {\it{isomorphic}} (in symbol, $L_1\cong L_2$), and the map $\phi:L_1\to L_2$ is an {\it{isomorphism}} if and only if $\phi$ is one-to-one and onto, and $a \leq b$ in $L_1$ if and only if $\phi(a) \leq \phi(b)$ in $L_2$. Algebraically, $\phi : L_1 \to L_2$ is an {\it{isomorphism}} if and only if $\phi$ is one-to-one and onto, and preserves both meet and join for any two elements. 

An element $b$ in $P$ {\it{covers}} an element $a$ in $P$ if $a<b$, and there is no element $c$ in $P$ such that $a<c<b$. Denote this fact by $a\prec b$, and say that pair $<a, b>$ is a $covering$ or an $edge$. If $a\prec b$ then $a$ is called a {\it{lower cover}} of $b$, and $b$ is called an {\it upper cover} of $a$. An element of a poset $P$ is called {\it{doubly irreducible}} if it has at most one lower cover and at most one upper cover in $P$. Let $Irr(P)$ denote the set of all doubly irreducible elements in the poset $P$. Let $Irr^*(P)= \{x \in Irr(P):x$ has exactly one upper cover and exactly one lower cover in $P\}$. The set of all coverings in $P$ is denoted by $E(P)$. The graph on a poset $P$ with edges as covering relations is called the {\it{cover graph}} and is denoted by $C(P)$. The number of coverings in a chain is called {\it{length}} of the chain.

The {\it nullity of a graph} $G$ is given by $m-n+c$, where $m$ is the number of edges in $G$, $n$ is the number of vertices in $G$, and $c$ is the number of connected components of $G$. Bhavale and Waphare \cite{bw} defined {\it nullity of a poset} $P$, denoted by $\eta(P)$, to be the nullity of its cover graph $C(P)$. For $a<b$, the interval $[a,b]=\{x\in P:a\leq x\leq b\}$, and $[a,b)=\{x\in P:a\leq x<b\}$; similarly, $(a,b)$ and ($a,b]$ can also be defined. For integer $n\geq 3$, {\it{crown}} is a partially ordered set $\{x_1,y_1,x_2,y_2,x_3,y_3,\ldots,x_n,y_n\}$ in which $x_i\leq y_i$, $y_i\geq x_{i+1}$, for $i=1,2,\ldots,n-1$, and $x_1\leq y_n$ are the only comparability relations. An element $x$ in a lattice $L$ is $join$-$reducible$($meet$-$reducible$) in $L$ if there exist $y,z \in L$ both distinct from $x$, such that $y\vee z=x (y\wedge z=x)$. An element $x$ in a lattice $L$ is $reducible$ if it is either join-reducible or meet-reducible. $x$ is $join$-$irreducible$($meet$-$irreducible$) if it is not join-reducible(meet-reducible); $x$ is $doubly$ $irreducible$ if it is both join-irreducible and meet-irreducible. The set of all doubly irreducible elements in $L$ is denoted by $Irr(L)$, and its complement in $L$ is denoted by $Red(L)$. An {\it{ear}} of a loopless connected graph $G$ is a subgraph of G which is a maximal path in which all internal vertices are of degree $2$ in $G$. {\it{Trivial~ear}} is an ear containing no internal vertices. A {\it{non-trivial ear}} is an ear which is not an edge. A vertex of a graph is called {\it pendant} if its degree is one.

\subsection{Dismantlable lattices and basic blocks}
In 1974, Rival \cite{ir} introduced and studied the class of dismantlable lattices. 
\begin{defn}\cite{ir}
A finite lattice of order $n$ is called {\it{dismantlable}} if there exists a chain $L_{1} \subset L_{2} \subset \ldots\subset L_{n}(=L)$ of sublattices of $L$ such that $|L_{i}| = i$, for all $i$.
\end{defn}

The following result is due to Kelly and Rival \cite{kr} (see Theorem 3.1 in \cite{kr}).
\begin{thm}\cite{kr}\label{crown}
A finite lattice is dismantlable if and only if it contains no crowns.
\end{thm}

The concept of \lq 1-sum\rq ~ is introduced by Bhavale and Waphare \cite{bw2}.
Let $P_1$ and $P_2$ be disjoint posets. Let $a \in P_1$. 
Define a partial order on $P = P_1 \cup P_2$ with respect to $a$ as follows.
For $x, y \in P$, we say that $x \leq y$ in $P$ if 
$x, y \in P_1$ and $x \leq y$ in $P_1$ or 
$x, y \in P_2$ and $x \leq y$ in $P_2$ or 
$x \in P_1, y \in P_2$ and $x \leq a$ in $P_1$.
It is easy to see that $P$ is a poset containing $P_1$ and $P_2$ as subposets. 
The procedure for obtaining $P$ in this way is called an {\it up 1-sum of $P_1$ with $P_2$ with respect to $a$}, 
denoted by $P_1 ]_a  P_2$. 
A diagram of $P$ is obtained by placing a diagram of $P_1$ and a diagram of $P_2$ side by side in such a way that the minimal elements of $P_2$ are at higher positions than $a$ and then by adding the coverings $<a, x>$ for all $x \in S$, the set of all minimal elements of $P_2$. This clearly gives $|E(P)| = |E(P_1)| + |E(P_2)| + |S|$.
Dually, one can define a {\it down 1-sum of two posets}. 
If $P$ is a down 1-sum of $P_1$ with $P_2$ with respect to $a$ in $P_1$ then we write $P = P_1 ]^a  P_2$. 
We call the element $a$ an {\it adjunct element} of the 1-sum. 
We say that $P$ is a {\it 1-sum of posets $P_1$ and $P_2$ with respect to an element $a \in P_1$} if $P$ is either an up 1-sum or a down 1-sum of $P_1$ and $P_2$ with respect to $a$.
A 1-sum $P_1 ]_a  P_2$ or $P_1 ]^a  P_2$ is called {\it trivial 1-sum} if $P_2$ is a chain and $a$ is respectively a maximal or a minimal element of $P_1$; otherwise, we say that the 1-sum is {\it non-trivial}.

The concept of {\it{adjunct operation of lattices}}, is introduced by Thakare, Pawar, and Waphare \cite{tpw}. 
Suppose $L_1$ and $L_2$ are two disjoint lattices and $(a, b)$ is a pair of elements in $L_1$ such that $a<b$ and $a\not\prec b$. Define the partial order $\leq$ on $L = L_1 \cup L_2$ with respect to the pair $(a,b)$ as follows: $x \leq y$ in $L$ if $x,y \in L_1$ and $x \leq y$ in $L_1$, or $x,y \in L_2$ and $x \leq y$ in $L_2$, or $x \in L_1,$ $ y \in L_2$ and $x \leq a$ in $L_1$, or $x \in L_2,$ $ y \in L_1$ and $b \leq y$ in $L_1$. 
It is easy to see that $L$ is a lattice containing $L_1$ and $L_2$ as sublattices. The procedure for obtaining $L$ in this way is called {\it{adjunct operation (or adjunct sum)}} of $L_1$ with $L_2$. We call the pair $(a,b)$ as an {\it{adjunct pair}} and $L$ as an {\it{adjunct}} of $L_1$ with $L_2$ with respect to the adjunct pair $(a,b)$ and write $L = L_1 ]^b_a L_2$. A diagram of $L$ is obtained by placing a diagram of $L_1$ and a diagram of $L_2$ side by side in such a way that the largest element $1$ of $L_2$ is at lower position than $b$ and the least element $0$ of $L_2$ is at the higher position than $a$ and then by adding the coverings $<1,b>$ and $<a,0>$. This clearly gives $|E(L)|=|E(L_1)|+|E(L_2)|+2$. A lattice $L$ is called an {\it{adjunct of lattices}} $L_1,L_2,\ldots,L_k,$ if it is of the form
$L = L_1 ]^{b_1}_{a_1} L_2 \cdots ]^{b_{k-1}}_{a_{k-1}} L_k$. 

The following result is due to Thakare, Pawar, and Waphare \cite{tpw} (see Theorem 2.2 in \cite{tpw}).
\begin{thm}\cite{tpw}\label{dac} 
A finite lattice is dismantlable if and only if it is an adjunct of chains.
\end{thm}

The following result is due to Bhavale and Aware \cite{ba} (see Corollary 2.5 in \cite{ba}).
\begin{cor} \cite{ba} \label{r-1}
A dismantlable lattice with $n$ elements has nullity $r-1$ if and only if it is an adjunct of $r$ chains.
\end{cor}

Thakare, Pawar, and Waphare \cite{tpw} defined a \textit{block} as a finite lattice in which the largest element is join-reducible and the least element is meet-reducible.

If $M$ and $N$ are two disjoint posets, the {\it{direct sum}} ({\it{see \cite{rp}}}), denoted by $M \oplus N$, is defined by taking the following order relation on 
$M\cup N$ : $x\leq y$ if and only if $x,y\in M$ and $x\leq y$ in $M$, or $x,y\in N$ and $x\leq y$ in $N$, or $x\in M,y\in N$. 
In general, $M \oplus N \neq N \oplus M$.
Also, if $M$ and $N$ are lattices then $|E(M\oplus N)|=|E(M)|+|E(N)|+1$.

\begin{rem}{\textnormal{Let $L$ be a finite lattice which is not a chain. Then $L$ contains a unique maximal sublattice which is a block, called the {\it{maximal block}}. The lattice $L$ has the form $C_1\oplus\textbf{B}$ or $\textbf{B}\oplus C_2$ or $L=C_1\oplus\textbf{B}\oplus C_2$, where $C_1,C_2$ are the chains and $\textbf{B}$ is the maximal block, hence $|E(L)|-|L|=|E(\textbf{B})|-|\textbf{B}|$}}.
\end{rem}

Bhavale and Waphare \cite{bw} introduced the following concepts namely, retractible element, basic retract, basic block, and basic block associated to a poset.
\begin{defn}\cite{bw}\label{rtrct} 
Let $P$ be a poset. Let $x \in Irr(P)$. Then $x$ is called a {\it{retractible}} element of $P$ if it satisfies either of the following conditions.
\begin{enumerate}
\item There are no $y,z \in Red(P)$ such that $y \prec x \prec z$.
\item There are $y,z \in Red(P)$ such that $y \prec x \prec z$ and there is no other directed path from $y$ to $z$ in $P$.
\end{enumerate}
\end{defn}

\begin{defn}\cite{bw} \label{brt}
A poset $P$ is a $basic$ $retract$ if no element of $Irr^{*}(P)$ is retractible in the poset $P$.
\end{defn}

\begin{defn}\cite{bw}\label{basicblock}
\textnormal{A poset $P$ is a {\it{basic block}} if it is one element or $Irr(P) = \emptyset$ or removal of any doubly irreducible element reduces nullity by one}.
\end{defn}

\begin{defn}\cite{bw}\label{bbas}
\textnormal{$B$ is a $basic$ $block$ $associated$ $to$ $a$ $poset$ $P$ if $B$ is obtained from the basic retract associated to $P$ by successive removal of all the pendant vertices.}
\end{defn}

The following result is due to Bhavale and Waphare \cite{bw} (see 3 of Theorem 3.5 in \cite{bw}).
\begin{thm}\cite{bw}\label{redb}
Let $B$ be a basic retract associated to a poset $P$. Then $Red(B) = Red(P)$.
\end{thm}

\subsection{Lattices in which reducible elements are comparable}\label{sec2}

Bhavale and Waphare \cite{bw} introduced the concept of RC-lattices as the class of lattices in which the reducible elements are all lying on a chain.
Using Theorem \ref{crown}, they have proved that RC-lattices are dismantlable.
Interestingly, the lattices on up to three reducible elements and having nullity up to two are RC-lattices.

\subsubsection{Counting of lattices on up to three reducible elements}

Although Thakare, Pawar, and Waphare \cite{tpw} counted up to isomorphism all lattices on two reducible elements, Bhavale and Aware \cite{ba} used a slightly different technique to count up to isomorphism all lattices on two reducible elements. Let $\mathscr{L}_r(n)$ denote the class of all non-isomorphic lattices $n$ elements such that every member of it contains exactly $r$ reducible elements. We denote the number of partitions of an integer $n$ into $k$ (non-decreasing and positive) parts by $P^k_n$. The following results are due to Bhavale and Aware \cite{ba} (see Theorem 3.5 and Theorem 4.16 in \cite{ba}), which provide the explicit formulae for the cases $r = 2$ and $r = 3$.

\begin{thm}\label{2r} \cite{ba}
For an integer $n \geq 4, \displaystyle |\mathscr{L}_2(n)| =\sum_{i=0}^{n-4}\sum_{k=0}^{n-i-4}(i+1)P^{k+2}_{n-i-2}$.
\end{thm}

\begin{thm} \cite{ba}
For an integer $n\geq 6$,\\
$\displaystyle |\mathscr{L}_3(n)|=\sum_{j=0}^{n-6}\sum_{k=1}^{n-j-5}\sum_{l=1}^{n-j-5}\sum_{i=1}^{n-j-l-4}2(j+1)P^{k+1}_{n-j-l-i-2}\\
~~~~~~~~~~~~+\sum_{j=0}^{n-6}\sum_{k=2}^{n-j-5}\sum_{r=5}^{n-j-2}\sum_{s=1}^{k-1}\sum_{i=1}^{r-4}2(j+1)P^{s+1}_{r-i-2}P^{k-s+1}_{n-j-r}\\
~~~~~~~~~~~~+\sum_{j=0}^{n-7}\sum_{k=1}^{n-j-6}\sum_{l=4}^{n-j-3}\sum_{t=1}^{k}(j+1)P^{t+1}_{l-2}P^{k-t+2}_{n-j-l-1}\\
~~~~~~~~~~~~+\sum_{j=0}^{n-8}\sum_{k=2}^{n-j-6}\sum_{r=1}^{n-j-7}\sum_{l=4}^{n-j-r-3}\sum_{t=1}^{k-1}(j+1) P^{t+1}_{l-2}P^{k-t+1}_{n-j-r-l-1}\\
~~~~~~~~~~~~+\sum_{j=0}^{n-8}\sum_{k=3}^{n-j-6}\sum_{r=2}^{n-j-7}\sum_{s=2}^{k-1}\sum_{l=4}^{n-j-r-3}\sum_{t=1}^{k-s}(j+1)P^{t+1}_{l-2}P^{k-s-t+2}_{n-j-r-l-1}P^{s}_{r}$.
\end{thm}

\subsubsection{Counting of RC-lattices containing up to $4$ reducible elements and having nullity up to $3$}

Let ${\mathscr L}_r^k (n)$ denote the class of all non-isomorphic RC-lattices on $n$ elements such that each lattice in it has nullity $k$ and contains $r$ reducible elements. 
Let ${\mathscr B}_r^k (m)$ denote the class of all non-isomorphic maximal blocks on $m$ elements such that each maximal block in it is an RC-lattice of nullity $k$ and contains $r$ reducible elements.
Note that, there is only one lattice, a chain, having nullity zero. Therefore ${\mathscr L}_0^0 (n)$ consists of the chain on $n$ elements. 

Bhavale and Waphare \cite{bw} proved that, if $L \in {\mathscr L}_r^k (n)$ then for fixed $k$, $2 \leq r \leq 2k$.
Therefore for $k=1$, $r=2$, and for $k=2$, $r=2,3,4$. Let ${\mathscr L}^2 (n) = {\mathscr L}_2^2 (n) \cup {\mathscr L}_3^2 (n) \cup {\mathscr L}_4^2 (n)$.
Let $<x>$ denote the nearest integer of a real number $x$.
The enumeration of all non-isomorphic lattices on $n$ elements and having nullity up to two was carried out by Thakare, Pawar, and Waphare \cite{tpw} (see Corollary 3.5 and Theorem 3.8 in \cite{tpw}). 

\begin{thm} \cite{tpw} \label{n1}
{\it For an integer} $n \geq 4$,
$$|{\mathscr L}_2^1 (n)| = \left \{ \begin{array}{lll} \frac{m(m-1)(4m+1)}{6} & \text{if} & n = 2m + 1; \\ & & \\ \frac{m(m-1) (4m-5)}{6} & \text{if} & n = 2m. \end{array} \right.$$
\end{thm}

\begin{thm}\cite{tpw} \label{n+1}
For an integer $n \geq 5$,
$|{\mathscr L}^2 (n)| =  \displaystyle\sum^{n-5}_{i=0} (i+1) |{\mathscr B}^2 (n-i)|,$ \\
where $ |{\mathscr B}^2 (j)| = \left \{ \begin{array}{lll}
< \frac{14{k^4} + 54{k^3} + 68{k^2} + 36 k + 9}{12} > & \text{if} & j = 2k+5; \\ & & \\
\lfloor \frac{(k+2)(7{k^3}+ 27{k^2} +31 k + 13 )}{6}\rfloor & \text{if} & j = 2k+6.\end{array}\right.$
\end{thm}

Recently, Aware and Bhavale \cite{AB} counted all non-isomorphic RC-lattices containing $4$ reducible elements and having nullity $3$.
The following result is due to Aware and Bhavale \cite{AB}.
\begin{thm} \cite{AB} \label{Bj43}
For $j\geq 7$, $|{\mathscr{B}}_4^3 (j)| = \displaystyle \sum_{p=1}^{j-6}\binom{j-p-2}{4}$\\
$+ \sum_{s=1}^{j-6}\sum_{r=1}^{j-s-5} \sum_{l=2}^{j-s-r-3}2(j-s-r-l-2)P_{l}^2$\\
$+ \sum_{t=1}^{j-7}\sum_{i=2}^{j-t-5}(i-1)P^{3}_{j-t-i-2}$\\
$+ \sum_{p=4}^{j-4}\sum_{t=1}^{j-p-3}tP_{j-p-t-1}^{2}P^2_{p-2}$\\
$+\sum_{t=1}^{j-7}\sum_{r=1}^{j-t-6}\sum_{l=1}^{j-t-r-5}\sum_{i=1}^{j-t-r-l-4}7P_{j-t-r-l-i-2}^{2}$\\
$+ \sum_{r=0}^{j-9}\sum_{p=5}^{j-r-4}2P_{p-2}^{3}P_{j-p-r-2}^{2}$\\
$+ \sum_{p=4}^{j-5}\sum_{l=1}^{j-p-4}\sum_{i=1}^{j-p-l-3}4P_{p-2}^{2}P_{j-p-l-i-1}^2$\\
$+ \sum_{r=1}^{j-8}\sum_{q=1}^{j-r-7}\sum_{l=4}^{j-q-r-3}2P_{l-2}^2P_{j-q-r-l-1}^2$\\
$+ \sum_{t=1}^{j-8}\sum_{m=0}^{j-t-8}\sum_{s=4}^{j-t-m-4}(j-t-m-7)P^{2}_{s-2}P^{2}_{j-t-m-s-2}$\\
$+ \sum_{p=7}^{j-3}\sum_{l=4}^{p-3}P_{j-p-1}^2P_{l-2}^2P_{p-l-1}^2$.
\end{thm}

The following result follows from Theorem \ref{Bj43}.
\begin{thm} \cite{AB} \label{Mainthm}
For $n\geq 7$, $|{\mathscr L}_4^3 (n)| = \displaystyle\sum_{i=0}^{n-7}(i+1)|{\mathscr{B}}_4^3 (n-i)|$.
\end{thm}

\noindent
In Section \ref{sec3}, we count the class of all non-isomorphic lattices on $n$ elements such that each lattice in it has nullity three, and at least two out of the four reducible elements in each lattice in it are incomparable.

\section{Lattices in which reducible elements are incomparable} \label{sec3}

\noindent
It is known that the reducible elements of a lattice of nullity up to two are all comparable (see \cite{tpw}). But a lattice of nullity at least three may not be RC-lattice. We define the class of RI-lattices as the class of lattices such that each lattice in it contains at least two incomparable reducible elements. 
By Theorem \ref{crown}, it follows that the lattices of nullity up to three, containing at most $7$ reducible elements are dismantlable, since cube ($2^3$) is the smallest lattice of nullity $3$, containing $8$ reducible elements, and which also contains the crown on $6$ reducible elements.
Therefore a lattice containing up to four reducible elements and having nullity up to three is dismantlable.

Let ${\mathcal L}_r^k (n)$ denote the class of all non-isomorphic dismantlable RI-lattices on $n$ elements such that each lattice in it contains $r$ reducible elements and has nullity $k$. 
Let ${\mathcal B}_r^k (m)$ denote the class of all non-isomorphic dismantlable maximal blocks on $m$ elements such that each maximal block in it is an RI-lattice of nullity $k$ and contains $r$ reducible elements.

Now in the following we prove that, there are three non-isomorphic basic blocks containing four reducible elements and having nullity three such that at least two reducible elements are incomparable.

\begin{prop} \label{i34}
If $B$ is the basic block associated to a lattice $L \in {\mathcal L}_4^3 (n)$ then $B \in \{B_1, B_2, B_3\}$ (see {\it Figure I}).
\end{prop}

\begin{center}
%TeXCAD Picture [Nullity-3-Fig9.pic]. Options:
%\grade{\on}
%\emlines{\off}
%\epic{\off}
%\beziermacro{\on}
%\reduce{\on}
%\snapping{\off}
%\pvinsert{% Your \input, \def, etc. here}
%\quality{8.000}
%\graddiff{0.005}
%\snapasp{1}
%\zoom{5.6569}
\unitlength 1mm % = 2.845pt
\linethickness{0.4pt}
\ifx\plotpoint\undefined\newsavebox{\plotpoint}\fi % GNUPLOT compatibility
\begin{picture}(96.783,31.437)(0,0)
\put(1.801,16.192){\circle{1.5}}
\put(1.801,23.44){\circle{1.5}}
\put(9.049,23.44){\circle{1.5}}
\put(16.297,8.944){\circle{1.5}}
\put(16.297,30.687){\circle{1.5}}
\put(23.545,23.44){\circle{1.5}}
\put(30.792,16.192){\circle{1.5}}
\put(30.792,23.44){\circle{1.5}}
\put(38.04,16.192){\circle{1.5}}
\put(38.04,23.44){\circle{1.5}}
\put(45.288,16.192){\circle{1.5}}
\put(52.536,8.944){\circle{1.5}}
\put(52.536,30.687){\circle{1.5}}
\put(59.784,16.192){\circle{1.5}}
\put(67.032,16.192){\circle{1.5}}
\put(67.032,23.44){\circle{1.5}}
%\emline(8.485,23.158)(1.945,16.44)
\multiput(8.485,23.158)(-.03371499,-.034626206){194}{\line(0,-1){.034626206}}
%\end
\put(16.617,30.229){\line(1,-1){6.364}}
\put(23.865,22.981){\line(1,-1){6.364}}
\put(1.591,22.981){\line(0,-1){6.01}}
\put(2.121,23.865){\line(2,1){13.435}}
%\emline(16.794,30.582)(30.052,23.688)
\multiput(16.794,30.582)(.064674108,-.033630536){205}{\line(1,0){.064674108}}
%\end
\put(16.794,9.192){\line(2,1){13.435}}
%\emline(2.121,15.733)(15.556,9.192)
\multiput(2.121,15.733)(.069252412,-.03371499){194}{\line(1,0){.069252412}}
%\end
\put(9.369,24.042){\line(1,1){6.364}}
\put(30.582,22.627){\line(0,-1){6.01}}
\put(38.36,15.733){\line(2,-1){13.435}}
%\emline(45.608,15.733)(51.972,9.192)
\multiput(45.608,15.733)(.033671599,-.034606921){189}{\line(0,-1){.034606921}}
%\end
%\emline(59.22,15.733)(53.033,9.546)
\multiput(59.22,15.733)(-.03362585,-.03362585){184}{\line(-1,0){.03362585}}
%\end
%\emline(66.468,15.733)(52.856,8.839)
\multiput(66.468,15.733)(-.066398751,-.033630536){205}{\line(-1,0){.066398751}}
%\end
%\emline(60.104,16.617)(66.468,23.158)
\multiput(60.104,16.617)(.033671599,.034606921){189}{\line(0,1){.034606921}}
%\end
%\emline(44.724,16.617)(38.184,23.334)
\multiput(44.724,16.617)(-.03371499,.034626206){194}{\line(0,1){.034626206}}
%\end
\put(37.83,22.804){\line(0,-1){6.187}}
\put(66.998,22.804){\line(0,-1){6.187}}
%\emline(38.184,23.865)(51.972,30.582)
\multiput(38.184,23.865)(.068942599,.03358742){200}{\line(1,0){.068942599}}
%\end
%\emline(53.033,30.582)(66.645,23.865)
\multiput(53.033,30.582)(.06805872,-.03358742){200}{\line(1,0){.06805872}}
%\end
\put(74.289,15.776){\circle{1.5}}
\put(74.289,23.024){\circle{1.5}}
\put(81.537,8.528){\circle{1.5}}
\put(81.537,23.024){\circle{1.5}}
\put(88.785,15.776){\circle{1.5}}
\put(88.785,30.272){\circle{1.5}}
\put(96.033,15.776){\circle{1.5}}
\put(96.033,23.024){\circle{1.5}}
%\emline(74.599,16.263)(81.14,22.804)
\multiput(74.599,16.263)(.03371499,.03371499){194}{\line(0,1){.03371499}}
%\end
%\emline(82.024,23.511)(88.388,30.052)
\multiput(82.024,23.511)(.033671599,.034606921){189}{\line(0,1){.034606921}}
%\end
%\emline(74.599,23.511)(88.034,30.229)
\multiput(74.599,23.511)(.06717484,.03358742){200}{\line(1,0){.06717484}}
%\end
\put(74.069,22.451){\line(0,-1){5.834}}
%\emline(74.599,15.203)(80.963,9.016)
\multiput(74.599,15.203)(.034586588,-.03362585){184}{\line(1,0){.034586588}}
%\end
%\emline(82.024,9.192)(88.388,15.556)
\multiput(82.024,9.192)(.033671599,.033671599){189}{\line(0,1){.033671599}}
%\end
%\emline(89.095,30.052)(95.636,23.511)
\multiput(89.095,30.052)(.03371499,-.03371499){194}{\line(1,0){.03371499}}
%\end
\put(96.166,22.451){\line(0,-1){5.834}}
%\emline(82.024,9.016)(95.459,15.733)
\multiput(82.024,9.016)(.06717484,.03358742){200}{\line(1,0){.06717484}}
%\end
%\emline(89.095,16.263)(95.459,22.627)
\multiput(89.095,16.263)(.033671599,.033671599){189}{\line(0,1){.033671599}}
%\end
\put(1.061,12.374){\makebox(0,0)[cc]{$a$}}
\put(30.582,12.198){\makebox(0,0)[cc]{$b$}}
\put(37.83,27.223){\makebox(0,0)[cc]{$a$}}
\put(67.175,27.4){\makebox(0,0)[cc]{$b$}}
\put(74.069,11.314){\makebox(0,0)[cc]{$a$}}
\put(95.989,27.047){\makebox(0,0)[cc]{$b$}}
\put(16.087,3.889){\makebox(0,0)[cc]{$B_1$}}
\put(52.502,3.889){\makebox(0,0)[cc]{$B_2$}}
\put(85.913,3.889){\makebox(0,0)[cc]{$B_3$}}
\end{picture}
\end{center}
\begin{center}
{\it Figure I}
\end{center}

\begin{proof}
Let$L \in {\mathcal L}_4^3 (n)$.  Suppose $B$ is the basic block associated to a lattice $L$.
Let $0, 1, a, b$ be the reducible elements of $L$. Note that by Theorem \ref{redb}, $Red(B) = Red(L)$. 
As at least two of them are incomparable, we have $a || b$. Also $a \wedge b = 0$ and $a \vee b = 1$. 
If $a$ (or $b$) is both meet as well as join reducible then $L$ is adjunct of at least $5$ chains, and hence nullity of $L$ is greater than or equal to $4$.
This is a contradiction. Therefore we are left with the following three possibilities.
\begin{enumerate}
\item If $a$ and $b$ both are meet reducible elements then $B$ is isomorphic to $B_1$ (see {\it Figure I}). 
\item If $a$ and $b$ both are join reducible elements then $B$ is isomorphic to $B_2$ (see {\it Figure I}). 
\item If without loss of generality, suppose $a$ is meet reducible element and $b$ is join reducible element, then $B$ is isomorphic to $B_3$ (see {\it Figure I}).
\end{enumerate}
\end{proof}

\section{Counting of RI-lattices containing $4$ reducible elements and having nullity $3$} 

For each $i, \; 1 \leq i \leq 3$, let ${\mathbb B}_i = \{{\bf B} \in {\mathcal L}_4^3 (n) ~:~ B_i$ is the basic block associated to ${\bf B} \}$.
Then ${\mathcal B}_4^3 (n) = {\mathbb B}_1 \dot\cup {\mathbb B}_2 \dot\cup {\mathbb B}_3$.
We define the class ${\mathcal L}'(n)$ as the subclass of ${\mathscr L}_2^1 (n)$, containing the lattices in which $1$ is the reducible element. 
Let ${\mathcal B}'(m)$ be the class of all maximal blocks in ${\mathcal L}'(n)$, where $m \leq n$.
In the following result, we obtain cardinality of the class ${\mathcal L}'(n)$.

\begin{lem}
For $n \geq 4$, $|{\mathcal L}'(n)| = \displaystyle\sum_{i=0}^{n-4} \lfloor \frac{n-i-2}{2} \rfloor$.
\end{lem}
\begin{proof}
Let $L \in {\mathcal L}'(n)$. Then $L = C \oplus {\bf B}$ where $C$ is a chain with $|C| = i \geq 0$ and ${\bf B} \in {\mathscr B}'(j)$ with $n=i+j$. Now $j \geq 4$. Therefore $i=n-j \leq n-4$. The proof follows from the fact that $|{\mathscr B}'(j)| = P_{j-2}^{2} = \lfloor \frac{j-2}{2} \rfloor$ for all $j \geq 4$.
\end{proof}

\begin{rem} \label{sj}
For any $j \geq 3$, let $S_j$ be the class of all non-isomorphic posets $Y$ on $j$ elements such that $Y = C ]_x C'$, where $C, C'$ are chains.
Then $Y \in S_j$ if and only if $Y \oplus \{1\} \in {\mathcal L}' (j+1)$. Therefore $|S_j| = |{\mathcal L}' (j+1)|$.
If $s_j = |S_j|$ for all $j$, then $s_j = |{\mathcal L}' (j+1)| = \displaystyle\sum_{i=0}^{j-3} \lfloor \frac{j-i-1}{2} \rfloor$.
\end{rem}

In the following result, we obtain cardinality of the class ${\mathbb B}_1$.

\begin{prop} \label{cb1}
For $n \geq 8$,\\
$|{\mathbb B}_1| = \left \{ \begin{array}{lll}
\displaystyle\sum_{n = i+j+2,~ i>j} s_i s_j, &~ \text{if} ~ n ~ is ~ odd; \\ 
\displaystyle\sum_{n = i+j+2,~ i>j} s_i s_j + \frac{s_{\frac{n-2}{2}} (s_{\frac{n-2}{2}}+1)}{2}, &~ \text{if} ~ n ~ is ~ even, \end{array}\right.$\\
where $s_i = \displaystyle\sum_{k=0}^{i-3} \lfloor \frac{i-k-1}{2} \rfloor$.
\end{prop}

\begin{proof}
Let ${\bf B} \in {\mathbb B}_1$. Then ${\bf B} \setminus \{0, 1\}$ is the disjoint union of two subposets, say $Y_1$ and $Y_2$ of ${\bf B}$, such that each one of them is an up 1-sum of two chains. By Theorem \ref{redb}, $Red({\bf B}) = Red(B_1)$. 
As $a, b \in Red(B_1)$, suppose $Y_1 = C_1 ]_a C_2$ and $Y_2 = C_3 ]_b C_4$ with $|Y_1| = i \geq 3$ and $|Y_2| = j \geq 3$, where $C_1, C_2, C_3$ and $C_4$ are chains. 
Suppose without loss of generality, ${\bf B} = (\{0\} \oplus Y_1 \oplus \{1\}) ]_0^1 Y_2$ with $|Y_1| = i \geq |Y_2| = j$ and $|{\bf B}| = n = i+j+2 \geq 8$.
It is clear that  $Y_1 \in S_i$ and $Y_2 \in S_j$. 
Let ${\bf B}' \in {\mathbb B}_1$ be such that ${\bf B}' = (\{0\} \oplus Y'_1 \oplus \{1\}) ]_0^1 Y'_2$. 
Then ${\bf B} \cong {\bf B}'$ if and only if $Y_1 \cong Y'_1$ and $Y_2 \cong Y'_2$.
Therefore, if $i > j$ then there are $\displaystyle\sum_{n = i+j+2} (|S_i| \times |S_j|)$ non-isomorphic maximal blocks in ${\mathbb B}_1$.
But if $i = j$ then $n$ must be even and it seems that there are $|S_i|^2$ maximal blocks (all may not be non-isomorphic).
In fact, there are $|S_i| \choose 2$ blocks which are counted twice, since $i=j$.
Therefore in the case when $i=j$, there are $|S_i|^2 - {{|S_i|} \choose 2} = \frac{|S_i|(|S_i|+1)}{2}$ non-isomorphic maximal blocks in ${\mathbb B}_1$. 
The proof follows from the fact that $s_i = |S_i| = \displaystyle\sum_{k=0}^{i-3} \lfloor \frac{i-k-1}{2} \rfloor$. 
\end{proof}

Using Corollary \ref{cb1}, in the following result we obtain cardinality of the class ${\mathbb B}_2$.

\begin{cor} \label{cb2}
For $n \geq 8$,\\
$|{\mathbb B}_2| = \left \{ \begin{array}{lll}
\displaystyle\sum_{n = i+j+2,~ i>j} s_i s_j, &~ \text{if $n$ is odd}; \\ 
\displaystyle\sum_{n = i+j+2,~ i>j} s_i s_j + \frac{s_{\frac{n-2}{2}} (s_{\frac{n-2}{2}}+1)}{2}, &~ \text{if $n$ is even}, \end{array}\right.$\\
where $s_i = \displaystyle\sum_{k=0}^{i-3} \lfloor \frac{i-k-1}{2} \rfloor.$
\end{cor}
\begin{proof}
Clearly $|{\mathbb B}_2| = |{\mathbb B}_1|$, since ${\bf B} \in {\mathbb B}_2$ if and only if the dual of ${\bf B}$, ${\bf B}^* \in {\mathbb B}_1$.
Thus the proof follows by Proposition \ref{cb1}.
\end{proof}

In the following result, we obtain cardinality of the class ${\mathbb B}_3$.
\begin{prop} \label{cb3}
For $n \geq 8$,
$|{\mathbb B}_3| = \displaystyle\sum_{n = i+j+2} s_i s_j$, where 
$s_i = \displaystyle\sum_{k=0}^{i-3} \lfloor \frac{i-k-1}{2} \rfloor$.
\end{prop}
\begin{proof}
Let ${\bf B} \in {\mathbb B}_3$. Then ${\bf B} \setminus \{0, 1\}$ is the disjoint union of two subposets, say $Y_1$ and $Y_2$ of ${\bf B}$ such that one of them is an up 1-sum of two chains, and the other is a down 1-sum of two chains. By Theorem \ref{redb}, $Red({\bf B}) = Red(B_3)$. 
As $a, b \in Red(B_3)$, suppose $Y_1 = C_1 ]_a C_2$ and $Y_2 = C_3 ]^b C_4$ with $|Y_1| = i \geq 3$ and $|Y_2| = j \geq 3$, where $C_1, C_2, C_3$ and $C_4$ are chains. 
Then either ${\bf B} = (\{0\} \oplus Y_1 \oplus \{1\}) ]_0^1 Y_2$ or 
${\bf B} = (\{0\} \oplus Y_2 \oplus \{1\}) ]_0^1 Y_1$ with 
$|{\bf B}| = n = i+j+2 \geq 8$. 
It is clear that $Y_1 \in S_i$, and the dual of $Y_2$, $Y_2^* \in S_j$.
Note that $Y_1 \oplus \{1\} \in {\mathscr L}^1 (i+1)$ and $(\{0\} \oplus Y_2)^* \in {\mathscr L}^1 (j+1)$.
Therefore, 
$|{\mathbb B}_3| = \displaystyle\sum_{n = i+j+2}(|S_i| \times |S_j|)$.
The proof follows from the fact that $s_k = |S_k| = \displaystyle\sum_{i=0}^{k-3} \lfloor \frac{k-i-1}{2} \rfloor$.
\end{proof}

In the following result, we obtain the number of all non-isomorphic maximal blocks on $n$ elements, having nullity three, and containing four reducible elements such that at least two of them are incomparable.

\begin{thm} \label{tn34}
For $n \geq 8$,\\
$|{\mathcal B}_4^3 (n)| = \left \{ \begin{array}{lll}
\displaystyle\sum_{n = i+j+2,~ i>j} 4 s_i s_j & \text{if $n$ is odd}; \\ 
\displaystyle\sum_{n = i+j+2,~ i>j} 4 s_i s_j + s_{\frac{n-2}{2}}
(2 s_{\frac{n-2}{2}}+1) & \text{if $n$ is even}, \end{array}\right.$\\
where $s_i = \displaystyle\sum_{k=0}^{i-3} \lfloor \frac{i-k-1}{2} \rfloor$.
\end{thm}
\begin{proof}
As $\{{\mathbb B}_i  : 1 \leq i \leq 3\}$ forms a partition of the class ${\mathcal B}_4^3 (n)$, we have 
$|{\mathcal B}_4^3 (n)| = |{\mathbb B}_1| + |{\mathbb B}_2| + 
|{\mathbb B}_3|$.
But by Corollary \ref{cb2},  $|{\mathbb B}_2| = |{\mathbb B}_1|$. 
Therefore $|{\mathcal B}_4^3 (n)| = 2 |{\mathbb B}_1| + |{\mathbb B}_3|$.
The remaining proof follows from Proposition \ref{cb1} and Proposition \ref{cb3}.
\end{proof}

Using Theorem \ref{tn34}, we have the following result.
\begin{thm}  \label{ri34}
For $n \geq 8$,
$|{\mathcal L}_4^3 (n)| = \displaystyle\sum_{i=0}^{n} |{\mathcal B}_4^3 (n-i)|$.
\end{thm}
\begin{proof}
Let $L \in {\mathcal L}_4^3 (n)$. Then $L=C \oplus {\bf B} \oplus C'$, where
${\bf B} \in {\mathcal B}_4^3 (j)$, and $C, C'$ are chains with $|C|+|C'|=i \geq 0$. Also $n=i+j \geq 8$. As $j \geq 8$, $i=n-j \leq n-8$. Further $i$ elements can be distributed on the chains $C$ and $C'$ in $i+1$ ways. Thus the proof follows from Theorem \ref{tn34}.
\end{proof}

\section*{Conclusion}
Using Theorem \ref{Mainthm} and Theorem \ref{ri34}, we obtain up to isomorphism all lattices on $n \geq 7$ elements, containing four reducible elements, and 
having nullity three.
We now raise the following two problems.\\
{\bf Problem 2.} Find up to isomorphism all lattices on $n$ elements, containing $r \geq 4$ comparable reducible elements, and having nullity $k \geq 4$.\\
{\bf Problem 3.} Find up to isomorphism all lattices on $n$ elements, containing $r \geq 4$ reducible elements, and having nullity $k \geq 4$.

\end{document}